\documentclass[12pt,a4paper]{article}

\usepackage{amsmath,amsfonts}
\usepackage{mathrsfs}

\newtheorem{theorem}{Theorem}
\newtheorem{lemma}{Lemma}

\def \Limsup{\mathop{\overline{\lim}}\limits}

\def\Pb{\mathbf{P}}

\def\Ex{\mathbf{E}}
\def\RR{\mathbb{R}}
\def\KK{\mathbb{K}}
\def\NN{\mathbb{N}}
\def\sgn{{\rm sgn}} 
\def\1{\mbox{1\hspace{-.25em}I}}
\begin{document}

\title{On Asymptotically Distribution  Free  Tests with Parametric Hypothesis
  for    Ergodic Diffusion Processes}
\author{M. \textsc{Kleptsyna},  Yu. A. \textsc{Kutoyants}\\
{\small Laboratoire de Statistique et Processus, Universit\'e du Maine}\\
{\small 72085 Le Mans, C\'edex 9, France}
}

\date{}

\maketitle
\begin{abstract}
We consider the problem of the construction of the asymptotically distribution
free test by the observations of ergodic diffusion process. It is supposedd
that under the basic hypothesis the trend coefficient depends on the finite
dimensional parameter and we study the Cram\'er-von Mises type statistics. The
underlying statistics depends on the deviation of the local time estimator
 from the invariant density with parameter replaced
by the maximum likelihood estimator. We propose a linear transformation which
yields the convergence of the test statistics to the integral of Wiener
process. Therefore the test based on this statistics is asymptotically
distribution free. 
\end{abstract}


\section{Introduction}

The goodness of fit (GoF) tests  occupy an important place in staistics because they
provide a bridge between the mathematical models and real data. Our work is
devoted to the problem of the construction of the GoF test in the case of
observation of ergodic diffusion process in the situation when the basic
hypothesis is composite parametric. We propose asymptotically distribution
free test, which is based on linear transformation of the normalized deviation
of the empirical density.

 Remind first the well-known properties of GoF tests in the statistics of
 i.i.d. observations $X_1,\ldots,X_n$. If we have to test the hypothesis
 ${\cal H}_0$ that their distribution function
 $F\left(x\right)=F_0\left(x\right)$ we can use (besides others) the
 Cram\'er-von Mises test $\hat\psi _n=\1_{\left\{\Delta _n>c_\varepsilon
   \right\}}$, where
$$
\Delta _n=n\int_{-\infty }^{\infty }\left[\hat
  F_n\left(x\right)-F_0\left(x\right)\right]^2 {\rm d}F_0\left(x\right),\qquad
\hat F_n\left(x\right)=\frac{1}{n}\sum_{j=1}^{n}\1_{\left\{X_j<x \right\}}.
$$ 
Remarcable property of this (and some other) tests is the fact that the 
statistics $\Delta _n$ under hypothesis ${\cal H}_0$ converges in distribution
$$
\Delta _n\Longrightarrow \Delta\equiv \int_{0}^{1}B\left(t\right)^2\,{\rm d}t,
$$
 where $B\left(t\right), 0\leq t\leq 1$ is a Brownien bridge. The tests with
the limit distribution not depending on the underlying model (here
$F_0\left(\cdot \right)$) are called {\it asymptotically distribution free }
(ADF). If we are interested in the construction of tests of asymptotically
fixed first type error $\varepsilon \in \left(0,1\right)$, i.e., the tests
$\bar\psi _n$ satisfying
$$
\lim_{n\rightarrow \infty }\Ex_{0 }\, \bar\psi _n=\varepsilon  ,
$$
 then for such tests the choice of the
threshold $c_\varepsilon $  can be done once for all problems with the same
limit distribution. Indeed, the threshold $c_\varepsilon $ for the test $\hat\psi _n$ is
solution of the equation  $ \Pb\left\{\Delta >c_\varepsilon\right\}=\varepsilon  
$, which is the same for all possible $F_0\left(\cdot \right)$.

If the basic hypothesis ${\cal H}_0$ is parametric:
$F\left(x\right)=F_0\left(\vartheta ,x\right)$, where $\vartheta \in \Theta
\subset \RR^d$ is an unknown parameter, then the situation changes and the limit 
distribution of the similar statistics
$$
\hat\Delta _n=n\int_{-\infty }^{\infty }\left[\hat
  F_n\left(x\right)-F_0\left(\hat\vartheta _n,x\right)\right]^2 {\rm
  d}F_0\left(\hat\vartheta _n,x\right)\Longrightarrow\hat \Delta , 
$$
($\hat\vartheta _n$ is the MLE) can be written in the following form
\begin{equation}
\label{bp}
\hat\Delta =\int_{0}^{1 }U\left(t\right)^2{\rm d}t,\qquad \quad
U\left(t\right)=B\left(t\right)-\left(\zeta 
  ,H\left(t\right)\right)
\end{equation}
where $\zeta=\zeta \left(\vartheta ,F_0\right) $ is a Gaussian vector and
$H\left(t \right)=H\left(\vartheta ,F_0,t\right)$ is some 
deterministic 
vector-function \cite{Dar}. If we decide to use  the test $\hat\psi _n=\1_{\left\{\hat\Delta
  _n>c_\varepsilon \right\}}$, then  we need to find such  $c_\varepsilon
=c_\varepsilon\left(\vartheta,F_0 \right) $ that $\Pb_\vartheta
\left(\hat\Delta >c_\varepsilon \right)=\varepsilon $,  verify  that $c_\varepsilon \left(\vartheta ,F_0 \right)$ is
continuous function of $\vartheta $ and  to put $\bar c_\varepsilon
=c_\varepsilon \left(\bar\vartheta_n ,F_0 \right)$, where $\bar\vartheta_n$ is
some consistent estimator of $\vartheta $ (say, MLE). Then it can be shown
that for the test $\hat \psi _n=\1_{\left\{\hat \Delta _n>\bar c_\varepsilon
 \right\}}$ we have
$$
\lim_{n\rightarrow \infty }\Ex_{\theta  }\, \hat\psi _n=\varepsilon \qquad
    {\rm for \quad all}\qquad \vartheta \in \Theta .
$$
We denote the  class of such tests as ${\cal K}_\varepsilon $. For a given
family $F_0\left(\cdot \right)$ the function $c_\varepsilon \left(\vartheta
,F_0 \right) $ can be found by numerical simulations. 
Of course, this problem becames much more complicate than the first one with
the simple basic hypothesis. More about GoF tests can be found, e.g.,
\cite{LR}, \cite{M10} or any other book on this subject.

Another possibility  is to find such 
transformormation $L\left[U_n\right]$ of the statistic  
$U_n\left(x\right)=\sqrt{n} \left(\hat F_n\left(x\right)-F(\hat \vartheta
_n,x)\right)$ that 
$$
\tilde\Delta _n=\int_{-\infty }^{\infty
}L\left[U_n\right]\left(x\right)^2{\rm d}F(\hat \vartheta 
_n,x)\Longrightarrow \tilde\Delta\equiv \int_{0}^{1}w_s^2\,{\rm d}s,\qquad
\Pb\left(\tilde\Delta >c_\varepsilon \right)=\varepsilon . 
$$
Then we will have the test $\tilde \psi _n=\1_{\left\{\tilde\Delta
  _n>c_\varepsilon \right\}}\in {\cal K}_\varepsilon$. Such linear
transformation was proposed in  \cite{Kh81}.

In our work we consider a similar problem of the construction of ADF GoF tests
by the observations of ergodic diffusion processes. We are given a stochastic
differential equation
\begin{equation}
\label{a}
{\rm d}X_s=S\left(X_s\right)\,{\rm d}s+\sigma \left(X_s\right)\,{\rm
  d}W_s,\quad X_0,\quad 0\leq s\leq T ,
\end{equation}
where $\sigma \left(x\right)^2>0$ is a known  function and we have to test a
composite basic hypothesis ${\cal H}_0$ that 
\begin{equation}
\label{b}
{\rm d}X_s=S\left(\vartheta ,X_s\right)\,{\rm d}s+\sigma \left(X_s\right)\,{\rm
  d}W_s,\quad X_0,\quad 0\leq s\leq T ,
\end{equation}
i.e., the trend coefficient is some known function $S\left(\vartheta
,x\right)$ which depends on the unknown parameter $\vartheta \in\Theta \subset
\RR^d$. Here and in the sequel we suppose that the initial value $X_0$ has the
distribution function of the invariant law of this ergodic diffusion process. 

Let us denote by $\hat F_T\left(x\right)$ and $\hat f_T\left(x\right)$
the {\it empirical distribution function} of the invariant law and the {\it empirical
density} (local time
estimator of the density $f\left(\vartheta ,x\right)$) defined by the relations
$$
\hat F_T\left(x\right)=\frac{1}{T}\int_{0}^{T}\1_{\left\{X_s<x\right\}}\,{\rm
  d}s,\qquad  \quad\hat f_T\left(x\right)=\frac{\Lambda
  _T\left(x\right)}{\sigma \left(x\right)^2T},  
$$ where $\Lambda _T\left(x\right)$ is the local time of the observed
diffusion process (see \cite{RY} for the definition and properties). Remind
that   we call  the random function $\hat f_T\left(x\right)$  {\it empirical
  density} because it is  the derivative of
empirical distribution function.

The Cram\'er-von Mises type statistics are based on $L_2$ deviations of these
estimators
\begin{align*}
\hat\eta _T\left(x\right)&=\sqrt{T}\left(\hat
  F_T\left(x\right)-F\left(\hat\vartheta _T,x\right)\right),\quad 
\hat\zeta _T\left(x\right)=\sqrt{T}\left(\hat
  f_T\left(x\right)-f\left(\hat\vartheta _T,x\right)\right),
\end{align*}
where $\hat\vartheta _T$ is the MLE of the parameter $\vartheta $. These
statistics  
can be introduced as follows
\begin{align*}
\hat\Delta _T=\int_{-\infty }^{\infty }\hat\eta _T\left(x\right) ^2{\rm
  d}F\left(\hat\vartheta _T,x\right),\qquad 
\hat\delta _T=\int_{-\infty }^{\infty }\hat\zeta _T\left(x\right) ^2{\rm
  d}F\left(\hat\vartheta _T,x\right).
 \end{align*}
 
Unfortunatelly,  the immediate  use of the tests $\hat\Psi _T=\1_{\left\{\hat\Delta
  _T>c_\varepsilon \right\}}$ and $\hat\psi _T=\1_{\left\{\hat\delta
  _T>d_\varepsilon \right\}}$ leads to the same problems as in the i.i.d. case,
i.e., the limit ($T\rightarrow \infty $)  distributions of these statistics
under hypothesis ${\cal H}_0$
depend  on the model $S\left(\cdot,\cdot  \right), \sigma \left(\cdot \right)$ and
on the true value $\vartheta$.

Moreover, despite the i.i.d. case, even if the basic hypothesis is simple
$\Theta =\left\{\vartheta _0\right\}$, the limit distributions depend on the
model defined by the functions $S\left(\vartheta_0 ,\cdot \right), \sigma
\left(\cdot \right)$. Therefore, even in this case of simple basic hypothesis
we have no ADF limits for these statistics.  This means that for each model we
have to find the threshold $c_\varepsilon $ separately. There are sevral ADF
GoF tests for the ergodic and ``small noise'' diffusion processes proposed,
for example, in the works \cite{DK},\cite{Kut11},\cite{NN}, but the links
between these tests and the ``traditional'' tests like Cram\'er-von Mises and
Kolmogorov-Smirnov (based on empirical distribution function) was not always
clear.

Recently in this problem (with simple hypothesis) there  was proposed a linear
transformation $L_1\left[\zeta _T\right]$  of the random function 
$$\zeta _T\left(x\right)=\sqrt{T} \left(\hat f_T\left(x\right)-f\left(\vartheta
_0,x\right)\right)$$ such that 
\begin{equation}
\label{c}
\delta _T=\int_{-\infty }^{\infty }\left[L_1\left[\zeta _T\right]\left(x\right)
  \right]^2{\rm d}F\left(\vartheta 
_0,x\right) \Longrightarrow  \int_{0}^{1}w\left(s\right)^2{\rm d}s 
\end{equation}
(see \cite{Kut12}). The proposed test statistics (after linear transformation
and some simplifications) is 
\begin{equation}
 \label{test1}
\tilde\delta _T=\int_{-\infty }^{\infty
}\left[\frac{1}{\sqrt{T}}\int_{0}^{T}\frac{\1_{\left\{X_s<x\right\}}}{\sigma
    \left(X_s\right)} \left[{\rm d}X_s-S\left(\vartheta _0,X_s\right){\rm d}s\right]
  \right]^2{\rm d}F\left(\vartheta 
_0,x\right)
\end{equation}
 with the same limit \eqref{c}. See as well \cite{NN}, where the similar
 statistics were used in the costruction of the Kolmogorov-Smirnov type  ADF test.

 Hence the test $\hat \psi _T=\1_{\left\{\tilde\delta
  _T>c_\varepsilon \right\}}$ is ADF (in the cas of simple basic hypothesis).

 The goal of this work is to present such 
linear transformation  $ L[\hat\zeta _T]$  of the random function
$\hat\zeta _T\left(x\right)$  that  
\begin{equation}
\label{cL}
\hat\delta _T=\int_{-\infty }^{\infty } L[\hat\zeta _T]\left(x\right)
  ^2{\rm d}F(\hat\vartheta
_T,x) \Longrightarrow  \int_{0}^{1}w\left(s\right)^2{\rm d}s .
\end{equation}

 Note that
the general case of ergodic diffusion process with shift (one-dimensional)
parameter was studied in \cite{NZ}. They showed that the limit
distribution of the Cram\'er-von Mises statistic does not depend on the unknown
(shift) parameter and therefore is {\it asymptotically parameter free}.

\section{Assumptions and Preliminaries}

We are given (under hypothesis ${\cal H}_0$) continuous time observations
$X^T=\left(X_s,0\leq s\leq T\right)$ of the diffusion process
\begin{equation}
\label{sde}
{\rm d}X_s=S\left(\vartheta ,X_s\right)\,{\rm d}s+\sigma \left(X_s\right)\,{\rm
  d}W_s,\quad X_0,\quad 0\leq s\leq T .
\end{equation}
We are going to study the GoF test based on the normalized difference 
\begin{align*}
&\sqrt{T}\left(\hat f_T\left(x\right)-f\left(\hat\vartheta
_T,x\right)\right)\\
&\qquad \qquad =\sqrt{T}\left(\hat
f_T\left(x\right)-f\left(\vartheta,x\right)\right)
-\left(\sqrt{T}\left(\hat\vartheta _T-\vartheta \right),\,\dot f\left(\vartheta
\right)\right)+ o\left(1\right) .
\end{align*}
 We need  three types of conditions. The first one (${\cal
  ES}$,${\cal RP}$ and ${\cal A}_0$)    provide the
existence of the solution of the equation \eqref{sde}, good ergodic properties of
the process $\left(X_s,s\geq 0\right)$ and allow to describe the asymptotic behavior of the
normalized difference $ \zeta _T\left(\vartheta,x\right)=\sqrt{T}\left(\hat
f_T\left(x\right)-f\left(\vartheta,x\right)\right)$. 

The regularity conditions ${\cal R}_1$ provide the properties of the MLE
$\hat\vartheta _T$ (consistency, asymptotic normality and stochastic
representation). The last condition ${\cal R}_2$ will help us to construct
the linear transformation $L \left[\cdot \right]$ of the process $\hat \zeta
_T\left(\cdot \right)$ to the Wiener process. Therefore, the test based on
this transformation is asymptotically distribution free.

We assume that the trend $S\left(\vartheta ,x\right),\vartheta \in
\Theta\subset \RR^d $ and diffusion $\sigma
\left(x\right)^2$ coefficients satisfy the following conditions.

${\cal ES}.$ {\it The function $S\left(\vartheta ,x\right),\vartheta
\in\Theta$ is locally bounded, the
function $\sigma \left(x \right)^2>0 $ is continuous  and for some $ C>0$
the condition
$$
x\,S\left(\vartheta ,x\right)+\sigma \left(x\right)^2\leq C\left(1+x^2\right)
$$
holds.}

By this condition the stochastic differential equation \eqref{sde} has a unique weak
solution  for all $\theta \in \Theta $  (see, e.g., \cite{Durrett}). 

${\cal RP}.$ {\it The functions $S\left(\vartheta ,\cdot \right)$  and $\sigma \left(x
\right)^2$ are such that for all $\vartheta \in \Theta $
$$
\int_{-\infty }^{x }\exp\left\{2\int_{0}^{x}\frac{S\left(\vartheta ,y\right)}{\sigma 
  \left(y\right)^2}{\rm d}y\right\} {\rm d}x\longrightarrow \pm \infty \quad
{as} \qquad x\longrightarrow \pm \infty 
$$
and}
$$
G\left(\vartheta \right)=\int_{-\infty }^{\infty }\sigma
\left(x\right)^{-2}\exp\left\{2\int_{0}^{x}\frac{S\left(\vartheta ,y\right)}{\sigma 
  \left(y\right)^2}{\rm d}y\right\} {\rm d}x<\infty .
$$

By  condition ${\cal RP}$  the diffusion process \eqref{sde} is recurrent positive (ergodic) with the
density of invariant law
$$
f\left(\vartheta ,x\right)=\frac{1}{G\left(\vartheta \right)\;\sigma
  \left(x\right)^2}\;\exp\left\{2\int_{0}^{x}\frac{S\left(\vartheta ,y\right)}{\sigma
  \left(y\right)^2}\; {\rm d}y\right\}.
$$
We suppose that the initial value $X_0$ has this density function, therefore
the observed process is stationary.

 Introduce the class ${\cal P}$ of  functions with polynomial
majorants
\begin{equation}
\label{p}
{\cal P}=\left\{h\left(\cdot \right):\quad \left|h\left(y\right)\right|\leq
C\left(1+\left|y\right|^p\right)\right\} .
\end{equation}
If the function $h\left(\cdot \right)$ depends on parameter $\vartheta $, then
we suppose that the constant $C$   in \eqref{p} does not depend on $\vartheta
$.  

The condition ${\cal RP} $ we strenghten by the following way.

${\cal A}_0.$ {\it The functions  $S\left(\vartheta ,\cdot \right), \sigma
  \left(\cdot \right)^{\pm 1} \in {\cal 
    P} $ and for all $\vartheta $}
$$
\Limsup_{\left|y\right|\rightarrow \infty 
}\;\sgn\left(y\right)\;\frac{S\left(\vartheta ,y\right)}{\sigma \left(y\right)^2} <0.
$$

The empirical distribution function $\hat F_T\left(x\right)$ and empirical
density $\hat f_T\left(x\right)$ by condition ${\cal A}_0$ are unbiased,
consistent, asymptotically normal and asymptotically efficient 
estimators of the functions $F\left(\vartheta ,x\right)$ and
$f\left(\vartheta ,x\right)$ respectively.  The random processes
\begin{align*} 
\eta _T\left(\vartheta,x\right)=\sqrt{T}\left(\hat F_T\left(x\right)-F\left(\vartheta ,
x\right)\right) ,\quad 
\zeta _T\left(\vartheta,x\right)=\sqrt{T}\left(\hat
f_T\left(x\right)-f\left(\vartheta , 
x\right)\right)
\end{align*}
converge to the Gaussian processes $\eta \left(\vartheta,x \right) $ and
$\zeta \left(\vartheta ,x \right) $, which admit the representations
\begin{align}
\label{3}
\eta \left(\vartheta,x \right)&=2\int_{-\infty }^{\infty
}\frac{F\left(\vartheta,y\right)F\left(\vartheta
  ,x\right)-F\left(\vartheta,y\wedge x\right)}{\sigma
  \left(y\right)\sqrt{f\left(\vartheta ,y\right)}}\;{\rm d}W\left(y\right),\\
\zeta \left(\vartheta ,x \right)&=2f\left(\vartheta ,x\right) \int_{-\infty }^{\infty
}\frac{F\left(\vartheta,y\right)-\1_{\left\{y>x\right\}}}{\sigma
  \left(y\right)\sqrt{f\left(\vartheta ,y\right)}}\;{\rm d}W\left(y\right).
\label{4}
\end{align}
Here $W\left(\cdot \right)$ is two-sided Wiener process. For the proofs see
\cite{Kut04}. These proofs 
 are based on the following representations  
\begin{align}
\label{df}
\eta _T\left(\vartheta,x\right)&=\frac{2}{\sqrt{T}}
\int_{0}^{T}\frac{F\left(\vartheta ,x\right)F\left(\vartheta ,X_s\right)
  -F\left(\vartheta ,x\wedge X_s\right)}{\sigma\left(y\right)
  \,f\left(\vartheta ,y\right)}\,{\rm d}W_s \nonumber \\
 &\qquad \quad
+\frac{2}{\sqrt{T}} \int_{X_0}^{X_T}\frac{F\left(\vartheta ,y\wedge
  x\right)-F\left(\vartheta ,y\right)F\left(\vartheta ,x\right)
}{\sigma\left(y\right)^2 \,f\left(\vartheta ,y\right)}\,{\rm d}y
\end{align}
and
\begin{align}
\label{lte}
\zeta _T\left(\vartheta,x\right)&=\frac{2f\left(\vartheta ,x\right)}{\sqrt{T}}
\int_{0}^{T}\frac{F\left(\vartheta ,X_s\right)-\1_{\left\{X_s>x\right\}}}{\sigma\left(y\right)
  \,f\left(\vartheta ,X_s\right)}\,{\rm d}W_s\nonumber\\
&\qquad \quad +\frac{2f\left(\vartheta ,x\right)}{\sqrt{T}}
\int_{X_0}^{X_T}\frac{\1_{\left\{y>x\right\}} -F\left(\vartheta ,y\right)}{\sigma\left(y\right)^2
  \,f\left(\vartheta ,y\right)}\,{\rm d}y.
\end{align}

It is easy to see that ${\cal A}_0$ implies ${\cal RP}$. Moreover, we can
verify that the condition ${\cal A}_0$ provides the equivalence of
the measures $\left\{\Pb_\vartheta ^{\left(T\right)},\vartheta \in \Theta
\right\}$ induced in the measurable space $\left({\cal
  C}\left[0,T\right],{\cal B}\right)$ of continuous on $\left[0,T\right]$
functions by the solutions of this equation with different $\vartheta
$ \cite{LS}. Hence, the likelihood ratio has the following form
\begin{align*}
L\left(\vartheta ,X^T\right)=\exp\left\{\int_{0}^{T}\frac{S\left(\vartheta
  ,X_s\right)}{\sigma \left(X_s\right)^2}\,{\rm
  d}X_s-\int_{0}^{T}\frac{S\left(\vartheta 
  ,X_s\right)^2}{2\,\sigma \left(X_s\right)^2}\,{\rm d}s \right\}
\end{align*}
and the MLE $\hat\vartheta _T$  is defined by the equation
$$
L\left(\hat\vartheta _T,X^T\right)=\sup_{\theta \in\Theta }L\left(\vartheta
,X^T\right) . 
$$ 

To study the tests we need to know the properties of the MLE $\hat\vartheta _T$ (in
the regular case). 

Below  and in the sequel the  dot means derivation w.r.t. $\vartheta $ and the
prime  means derivation
w.r.t. $x$, i.e.; $\dot S\left(\vartheta
,x\right)$ is $d$-vector and  $ \ddot S \left(\vartheta
,x\right) $ is a $d\times d$ matrix.
 The information matrix is
$$
{\rm I}\left(\vartheta \right)=\Ex_{\vartheta} \left(\frac{\dot S\left(\vartheta,\xi\right)
  \;\dot S\left(\vartheta,\xi\right)^*  }{\sigma \left(\xi \right)^2}\right),
$$
where * means transposition and $\xi $ is the r.v. with the invariant density
function $f\left(\vartheta ,x\right)$. The scalar product in $\RR^d$ we denote
as 
$\langle\cdot ,\cdot\rangle $. 

We have two types of {\it Regularity conditions}. 

${\cal R}_1.$
\begin{itemize}
\item  {\it The set  $\Theta $ is an  open and  bounded subset of $\RR^d$.}
 \item  {\it The function $S\left(\vartheta ,x\right)$    has continuous
   derivatives w.r.t. $\vartheta $ such that }
$$
\dot S\left(\vartheta
,x\right),\;  \ddot S \left(\vartheta
,x\right)\in {\cal P}.
$$
\item {\it The information matrix is uniformly nondegerate 
$$ 
\inf_{\vartheta \in \Theta }\inf_{\left|\lambda \right|=1,\lambda \in
    \RR^d} \lambda ^*{\rm I}\left(\vartheta \right)\lambda >0
$$ 
and for any compact  $\KK\subset \Theta $, any $\vartheta _0\in \Theta $ and any $\nu >0$ 
$$
\inf_{\vartheta \in \KK} \inf_{\left|\vartheta -\vartheta _0\right|>\nu
}\Ex_{\vartheta_0} \left(\frac{S\left(\vartheta,\xi\right)
  -S\left(\vartheta_0,\xi\right)  }{\sigma \left(\xi \right)}\right)^2>0. 
$$
 }
\end{itemize}
Here $\xi $ is a random variable with the density function $f\left(\vartheta
_0,x\right)$.
By the conditions ${\cal A}_0$  and ${\cal R}_1$ the MLE is consistent,
asymptotically normal
$$
\sqrt{T}\left(\hat\vartheta _T-\vartheta \right)\Longrightarrow {\cal
  N}\left(0, {\rm I}\left(\vartheta \right)^{-1}\right),
$$
 we have the convergence of all polynomial moments and this estimator is
 asymptotically efficient (see \cite{Kut04} for details). Moreover, the MLE
 admits the representation 
\begin{equation}
\label{mlE}
\sqrt{T}\left(\hat\vartheta _T-\vartheta \right)=\frac{{\rm I}\left(\vartheta
\right)^{-1}}{\sqrt{T}}\int_{0}^{T} \frac{\dot S\left(\vartheta
  ,X_s\right)}{\sigma \left(X_s\right)}\;{\rm
  d}W_s\,\left(1+o\left(1\right)\right) .
\end{equation}

Let us introduce the matrix
$$
N\left(\vartheta ,y\right)={\rm I}\left(\vartheta
\right)^{-1}\int_{y}^{\infty } \frac{\dot S\left(\vartheta ,z\right)\; \dot
  S\left(\vartheta ,z\right)^*}{\sigma \left(z\right)^2} \,f\left(\vartheta
,z\right)\,{\rm d}z  .
$$
Note that $N\left(\vartheta ,-\infty \right)=I_d $, where $I_d$ is the unit $d\times d$ matrix.

The next regularity condition is \\
${\cal R}_2.$ 
\begin{itemize}
\item {\it The functions $\dot S\left(\vartheta ,x\right)$ and $\sigma
  \left(x\right)$ have continuous derivatives w.r.t. $x$ 
$$ \dot S' \left(\vartheta
,x\right), \;\sigma '\left(x\right) \quad \in \quad {\cal P}.
$$
\item The matrix $N\left(\vartheta ,y\right) $ for any $y$ is uniformly non degenerate}
$$
\inf_{\vartheta \in \Theta }\inf_{\left|\lambda \right|=1,\lambda \in \RR^d}
\lambda ^*N\left(\vartheta ,y\right)\lambda >0. 
$$
\end{itemize}
Let us remind what happens in the case of simple basic hypothesis, say,
$\vartheta =\vartheta _0$.  Using the representation \eqref{df} and
\eqref{lte} it is shown that  the 
corresponding statistics have the following limits
\begin{align*}
\Delta _T&=T\int_{}^{ }\left[\hat F_T\left(x\right)-F\left(\vartheta
  _0,x\right)\right] ^2{\rm d}F\left(\vartheta _0,x\right) \Longrightarrow
\int_{}^{ }\eta \left(\vartheta _0,x\right)^2{\rm d}F\left(\vartheta
_0,x\right) ,\\ 
\delta _T&=T\int_{}^{ }\left[\hat
  f_T\left(x\right)-f\left(\vartheta _0,x\right)\right] ^2{\rm
  d}F\left(\vartheta _0,x\right)\Longrightarrow \int_{ }^{}\zeta
\left(\vartheta _0,x\right)^2{\rm d}F\left(\vartheta _0,x\right).
 \end{align*}
 Therefore the tests
based on these two statistics are not ADF. To construct the ADF test we put
$$
\mu _0\left(\vartheta _0,x\right) =\frac{\zeta \left(\vartheta
  _0,x\right)}{2f\left(\vartheta_0 ,y\right)}=\int_{-\infty }^{\infty
}\frac{F\left(\vartheta_0 ,y\right)-\1_{\left\{y>x\right\}} }{\sigma
  \left(y\right)\sqrt{f\left(\vartheta_0 ,y\right)}} \;{\rm
  d}W\left(y\right),
$$
and note that by the CLT
$$ \frac{\zeta _T\left(\vartheta_0 ,y\right)}{2f\left(\vartheta_0
  ,y\right)}\Longrightarrow \mu _0\left(\vartheta _0,x\right).
$$ 
 Further, we have the convergence
\begin{align}
L_1\left[\zeta _T\left(\vartheta_0  \right)\right]\left(x\right)&=\int_{-\infty
}^{x}  \sigma \left(y\right)f\left(\vartheta_0 ,y\right){\rm d}\left[\frac{\zeta
    _T\left(\vartheta_0 ,y\right)}{2f\left(\vartheta_0 ,y\right)}\right]\nonumber\\
&
=\frac{1}{\sqrt{T}}\int_{0}^{T}\1_{\left\{X_s<x\right\}}\,{\rm d}W_s+o\left(1\right)\Longrightarrow 
w\left({F\left(\vartheta_0 ,x\right)}\right) .
\label{L1}
\end{align}
Hence
\begin{align*}
\bar\delta _T=&\int_{-\infty }^{\infty }L_1\left[\zeta _T\left(\vartheta_0
  \right)\right]\left(x\right)^2{\rm d}F\left(\vartheta_0
,x\right)\\
&\quad 
\Longrightarrow \int_{-\infty }^{\infty }w\left({F\left(\vartheta_0 ,x\right)}\right)^2\,{\rm
  d}F\left(\vartheta_0 ,x\right)=\int_{0}^{1}w\left(s\right)^2\,{\rm
  d}s
\end{align*}
and the test $\bar \psi _T=\1_{\left\{ \bar\delta _T>c_\varepsilon \right\}}$
is ADF (see the details in \cite{Kut12}). 

Moreover,  we can define an asymptotically equivalent
test $\tilde\psi_T=\1_{\left\{\tilde \delta _T>c_\varepsilon \right\}} $, where
\begin{equation}
\label{h0}
\tilde \delta _T=\int_{-\infty }^{\infty }\left[
  \frac{1}{\sqrt{T}}\int_{0}^{T}\frac{\1_{\left\{X_s<x\right\}}}{\sigma
    \left(X_s\right)}\;\left[{\rm d}X_s-S\left(\vartheta_0 ,X_s\right)\,{\rm
      d}s\right]\right]^2 {\rm d}F\left(\vartheta_0 ,x\right)
\end{equation}
and this test as well is ADF.

\section{Main Result}
Remind that the value of parameter $\vartheta $ is unknown that is why  we replace
$\vartheta $ by its MLE $\hat\vartheta _T$ and
 our goal is to find the  transformations 
$$
L\left[\eta
  _T\left(\hat\vartheta _T,\cdot \right)\right]\left(x\right),\qquad L\left[\zeta
  _T\left(\hat\vartheta _T,\cdot \right)\right]\left(x\right)
$$ of the statistics $\eta _T(\hat\vartheta _T,x
)=\sqrt{T}\left(\hat F_T\left(x\right)-F(\hat\vartheta
_T,x)\right)$ and $\zeta _T\left(\hat\vartheta _T,x
\right)=$ \\
$\sqrt{T}\left(\hat f_T\left(x\right)-f\left(\hat\vartheta
_T,x\right)\right)$ such that the GoF tests constructed on it will be ADF. First note that we have
equality
$$
\left[\eta  _T(\hat\vartheta _T,x)\right]'=\zeta
_T(\hat\vartheta _T,x ) ,
$$ 
therefore if we  find
this  transformation for $\zeta
_T(\hat\vartheta _T,\cdot  ) $ then we obtain it for $\eta
_T(\hat\vartheta _T,\cdot ) $ too.

Moreover, we show that the linear transformation \eqref{L1} of 
$$
\mu  _T\left(\hat\vartheta _T,x \right)=\frac{\sqrt{T}(\hat
f_T\left(x\right)-f(\hat\vartheta _T,x))}{2f(\hat\vartheta _T,x) },\qquad x\in
\RR 
$$ 
 gives us statistic which is asymptotically equivalent to the  statistic
$$
\xi _T\left(\hat\vartheta
_T,x\right)=\frac{1}{\sqrt{T}}\int_{0}^{T}\frac{\1_{\left\{X_s<x\right\}}}{\sigma  
  \left(X_s\right)} \;\left[{\rm d}X_s-S(\hat\vartheta _T,X_s){\rm
    d}s\right].
$$
Therefore our  ADF test will be based on the statistic $\xi _T(\hat\vartheta
_T,x) $, which is much easier to calculate.

Introduce the random vector
\begin{equation}
\label{ka}
\Delta \left(\vartheta \right)=\int_{-\infty }^{\infty }\frac{\dot
  S\left(\vartheta ,y\right)}{\sigma \left(y\right)}\sqrt{f\left(\vartheta
  ,y\right)}\,{\rm d}W\left(y\right)\quad \sim \quad {\cal N}\left(0, {\rm
  I}\left(\vartheta \right)\right) 
\end{equation}
and the Gaussian function 
$$
\mu \left(\vartheta ,x\right)=  \mu _0 \left(\vartheta ,x\right) -2^{-1}\langle{\rm
  I}\left(\vartheta 
\right)^{-1} \Delta \left(\vartheta \right),\frac{\partial \ell\left(\vartheta
  ,x\right)}{\partial \theta }\rangle ,\qquad x\in \RR,
$$
where $\ell\left(\vartheta
  ,x\right)=\ln f\left(\vartheta ,x\right)$
and $\langle\cdot ,\cdot \rangle$ is the scalar product in $\RR^d$. Further, let
us put $s=F\left(\vartheta ,y\right)$, $t=F\left(\vartheta ,x\right)$,  define
the vector function 
$$
h\left(\vartheta,s\right)={\rm I}\left(\vartheta
\right)^{-1/2} \frac{\dot S \left(\vartheta,F^{-1}\left(\vartheta ,s\right)
\right)}{\sigma \left(F^{-1}\left(\vartheta ,s\right)\right)},\qquad
\int_{0}^{1}h\left(\vartheta ,s\right)^*h\left(\vartheta ,s\right){\rm 
  d}s =1,
$$
and Gaussian process 
\begin{equation}
\label{U}
U\left(t\right)=w\left(t\right)-\langle \int_{0}^{1}h\left(\vartheta ,s\right){\rm
  d}w\left(s\right), \int_{0}^{t}h\left(\vartheta ,s\right){\rm
  d}s\rangle,
\end{equation}
where $w\left(s\right), 0\leq s\leq 1$ is some Wiener process. 
Here $F^{-1}\left(\vartheta ,s\right) $ is the function inverse to
$F\left(\vartheta ,y\right)$, i.e., the solution $y$ of the equation
$F\left(\vartheta ,y\right)=s $. Below
$u\left(x\right)=U\left(F\left(\vartheta ,x\right)\right)$. 

\begin{theorem}
\label{T1}
Let the conditions ${\cal ES}, {\cal A}_0$ and ${\cal R}_1$ be fulfilled, then
\begin{align}
\label{f-d}
\mu  _T\left(\hat\vartheta _T,x \right) \Longrightarrow \mu \left(\vartheta
,x\right),\qquad \xi _T\left(\hat\vartheta _T,x \right) \Longrightarrow u\left(x\right)
, 
\end{align}
and
\begin{align}
\label{repr}
\int_{-\infty }^{x}\sigma \left(y\right)f\left(\vartheta ,y\right){\rm d}\mu
\left(\vartheta ,y\right)=u\left(x\right). 
\end{align}
\end{theorem}
{\bf Proof.} Using the consisteny of the MLE we can write 
\begin{align*}
\zeta _T\left(\hat\vartheta _T,x\right)&=\sqrt{T}\left(\hat
f_T\left(x\right)-f\left(\vartheta ,x\right)
\right)+\sqrt{T}\left(f\left(\vartheta ,x\right)-f(\hat\vartheta_T
,x) \right)\\
&=\zeta _T\left(\vartheta ,x\right)-\langle\sqrt{T}(\hat\vartheta_T-\vartheta
) ,\frac{\partial f\left(\vartheta ,x\right)}{\partial \vartheta }
\rangle +o\left(1\right)
. 
\end{align*}

 The slight modification of the proof of the Theorem 2.8 in \cite{Kut04}
 allows us to verify the joint asymptotic normality of  $\zeta
 _T\left(\vartheta ,x\right)$ and $\sqrt{T}\left(\hat\vartheta _T-\vartheta \right) $ as follows.
 Let us denote  $\Delta _T\left(\vartheta ,X^T\right) $  the vector score function
$$
\Delta _T\left(\vartheta ,X^T\right)=\frac{1}{\sqrt{T}}\int_{0}^{T}\frac{\dot
  S\left(\vartheta ,X_s\right)}{\sigma \left(X_s\right)}\;{\rm d}W_s. 
$$ 
The behavior of the MLE is described in \cite{Kut04} through the weak convergence of the
normalized likelihood ratio
$$
 Z_T\left(u\right)\equiv
\frac{L\left(\vartheta+\frac{u}{\sqrt{T}},X^T\right)}{L\left(\vartheta,X^T\right)}=\exp\left\{\langle  
u,\Delta
_T\left(\vartheta ,X^T\right)\rangle 
-\frac{1}{2} u^*{\rm I}\left(\vartheta
\right)u+o\left(1\right)\right\}.
$$
 By the central limit theorem for stochastic integrals we have the joint asymptotic
 normality: for any $\left(\lambda ,\nu \right) \in \RR^{1+d}$
$$
\lambda\, \zeta _T\left(\vartheta ,x \right)+\langle\nu, \Delta
_T\left(\vartheta ,X^T\right)\rangle\Longrightarrow 
\lambda\,\zeta\left(\vartheta ,x \right)+\langle\nu ,\Delta \left(\vartheta
\right) \rangle  .
$$
Hence following the proof of the mentioned above Theorem 2.8 we obtain the joint
convergence 
$$
\left(\zeta _T\left(\vartheta ,x \right), Z_T\left(\cdot
\right)\right)\Longrightarrow \left(\zeta_0\left(\vartheta ,x \right),
Z\left(\cdot \right)\right) ,
$$
where
$$
Z\left(u\right)=\exp\left\{\langle 
u,\Delta
\left(\vartheta \right)\rangle 
-\frac{1}{2} u^*{\rm I}\left(\vartheta
\right)u\right\},\qquad u\in \RR^d.
$$
This joint convergence yields the joint asymptotic normality 
$$ \left(\zeta _T\left(\vartheta ,x \right), \sqrt{T}(\hat\vartheta
_T-\vartheta ) \right)\Longrightarrow \left(\zeta \left(\vartheta ,x 
\right),{\rm I}\left(\vartheta \right)^{-1}\Delta \left(\vartheta \right)\right)
$$
with the same  Wiener process $W\left(\cdot \right)$ in \eqref{4} and
\eqref{ka}.

Now the convergence \eqref{f-d} follows from the consisteny of the MLE,
because $f(\hat\vartheta _T,x)\rightarrow f\left(\vartheta ,x\right)$.

Therefore the limit $\mu \left(\vartheta ,x\right) $ of $\mu _T\left(\vartheta
,x\right)$ can be written as
\begin{align*}
\int_{-\infty }^{\infty
}\left[\frac{F\left(\vartheta
    ,y\right)-\1_{\left\{y>x\right\}}-\langle\left[2{\rm I}\left(\vartheta 
\right)\right]^{-1} \dot S \left(\vartheta,y \right), \dot \ell\left(\vartheta
      ,x\right)\rangle f\left(\vartheta ,y\right)}{\sigma 
  \left(y\right)\sqrt{f\left(\vartheta ,y\right)}} \right]{\rm d}W\left(y\right) .
\end{align*}

Let us consider the linear transformation of $\mu \left(\vartheta ,\cdot \right)$
following \eqref{L1}:
\begin{align*}
L_1\left[\mu \right]\left(x \right)=\int_{-\infty }^{x}\sigma
\left(y\right)f\left(\vartheta ,y\right)\,{\rm d}\mu \left(\vartheta ,y\right) .
\end{align*}
Remind the details of this transformation  from \cite{Kut12}.  Denote
\begin{align*}
F\left(\vartheta ,y\right)&=s,\quad
a\left(\vartheta ,s\right)=\sigma \left(F^{-1}\left(\vartheta ,s\right)
\right),\quad 
b\left(\vartheta ,s\right)= f \left(\vartheta ,F^{-1}\left(\vartheta ,s\right) \right).
\end{align*}
Then we can write 
\begin{align*}
&\int_{-\infty }^{\infty
}\frac{F\left(\vartheta ,y\right)-\1_{\left\{y>x\right\}} }{\sigma
  \left(y\right)\sqrt{f\left(\vartheta ,y\right)}} \;{\rm
  d}W\left(y\right)\\
&\qquad  =\int_{-\infty }^{\infty 
}\frac{\left[F\left(\vartheta ,y\right)-\1_{\left\{F\left(\vartheta
      ,y\right)>F\left(\vartheta ,x\right)\right\}}\right] }{\sigma 
  \left(y\right){f\left(\vartheta ,y\right)}} \sqrt{f\left(\vartheta
  ,y\right)}\;{\rm d}W\left(y\right)\\
&\qquad  =\int_{0 }^{1 
}\frac{\left[s-\1_{\left\{s>t\right\}}\right]
}{a\left(\vartheta ,s\right)b\left(\vartheta ,s\right)} \;{\rm d}w\left(s\right)  \\
&\qquad =\int_{0 }^{t
}\frac{s
}{a\left(\vartheta ,s\right)b\left(\vartheta ,s\right)} \;{\rm d}w\left(s\right)+\int_{t }^{1 
}\frac{s-1
}{a\left(\vartheta ,s\right)b\left(\vartheta ,s\right)} \;{\rm d}w\left(s\right) \\
&\qquad=v\left(\vartheta ,t\right),\qquad 0<t<1,
\end{align*}
where  $w\left(s\right),0\leq s\leq 1$ is the following  Wiener process 
$$
w\left(s \right)=\int_{-\infty
}^{F^{-1}\left(\vartheta ,s\right)}\sqrt{f\left(\vartheta ,y\right)}\,{\rm d}W\left(y\right). 
$$
 Note that $v\left(\vartheta ,0\right)=\infty $ ($x=-\infty )$ and
$v\left(\vartheta ,1\right)=\infty $ ($x=+\infty )$. Therefore we define this
differential and the corresponding integrals below for $t\in \left(\nu ,1-\nu
\right)$ with small $\nu >0$ and in the sequel $\nu \rightarrow 0$
($x\rightarrow \pm \infty $).

Hence
$$
{\rm d}\mu_0 \left(\vartheta ,y\right)={\rm d}v\left(\vartheta ,s\right)=\frac{1
}{a\left(\vartheta ,s\right)b\left(\vartheta ,s\right)} \;{\rm d}w\left(s\right)
$$
and
$$
\int_{-\infty }^{x}\sigma
\left(y\right)f\left(\vartheta ,y\right)\,{\rm d}\mu_0 \left(\vartheta ,y\right)= \int_{0}^{t}a\left(\vartheta ,s\right)b\left(\vartheta ,s\right)\,{\rm
  d}v\left(\vartheta ,s\right)=w\left(t\right) .
$$
To calculate the second term note that
$$
\dot \ell \left(\vartheta ,x\right)=-\frac{\dot G\left(\vartheta
  \right)}{G\left(\vartheta \right)} +2\int_{0}^{x}\frac{\dot S\left(\vartheta
  ,y\right)}{\sigma \left(y\right)^2}{\rm d}y 
.$$
Therefore
\begin{align*}
\int_{-\infty }^{x}\sigma \left(y\right)f\left(\vartheta ,y\right){\rm d}\dot
\ell \left(\vartheta ,y\right) =2\int_{-\infty }^{x} \frac{\dot
  S\left(\vartheta ,y\right)}{\sigma \left(y\right)}  f\left(\vartheta
,y\right){\rm d}y 
\end{align*}
and
\begin{align*}
&\int_{-\infty }^{x}\sigma \left(y\right)f\left(\vartheta ,y\right){\rm d}\mu
\left(\vartheta ,y\right)=w\left(F\left(\vartheta ,x\right)\right)\\
&\qquad -\langle {\rm
  I}\left(\vartheta \right)^{-1/2} \int_{-\infty }^{\infty }\frac{\dot 
  S\left(\vartheta ,y\right)}{\sigma \left(y\right)}{\rm
  d}w\left(F\left(\vartheta ,y\right)\right), {\rm 
  I}\left(\vartheta \right)^{-1/2} \int_{-\infty }^{x }\frac{\dot
  S\left(\vartheta ,y\right)}{\sigma \left(y\right)}{\rm d}F\left(\vartheta
,y\right)\rangle\\
&\quad =U\left(F\left(\vartheta
  ,x\right)\right)=w\left(t\right)-\langle \int_{0}^{1}h\left(\vartheta ,s\right){\rm
  d}w\left(s\right), \int_{0}^{t}h\left(\vartheta ,s\right){\rm
  d}s\rangle.
\end{align*}
Further, we have 
\begin{align}
&\xi _T\left(\hat\vartheta _T,x\right)=\frac{1}{\sqrt{T}}\int_{0}^{T}\frac{\1_{\left\{X_s<x\right\}}}{\sigma
  \left(X_s\right)} \;\left[{\rm d}X_s-S(\vartheta,X_s){\rm
    d}s\right]\nonumber\\
& \qquad +\frac{1}{\sqrt{T}}\int_{0}^{T}\frac{\1_{\left\{X_s<x\right\}}}{\sigma
  \left(X_s\right)} \;\left[S(\vartheta,X_s)-S(\hat\vartheta _T,X_s)\right]{\rm
    d}s\nonumber\\
&\quad =\frac{1}{\sqrt{T}}\int_{0}^{T}\1_{\left\{X_s<x\right\}}{\rm
  d}W_s-\langle\left(\hat\vartheta _T- \vartheta
  \right),\int_{0}^{T}\frac{\1_{\left\{X_s<x\right\}}\dot S(\vartheta,X_s)}{\sqrt{T}\sigma 
  \left(X_s\right)} \;{\rm
    d}s\rangle+o\left(1\right)\nonumber\\
&\quad \Longrightarrow w\left(F\left(\vartheta ,x\right)\right)-\langle{\rm
    I}\left(\vartheta \right)^{-1}\Delta \left(\vartheta \right),\int_{-\infty
  }^{x}\frac{\dot S\left(\vartheta ,y\right)}{\sigma \left(y\right)}\; {\rm
    d}F\left(\vartheta ,y\right)\rangle=u\left(x\right).
\label{55}
\end{align}

 It can be shown that 
$$
L_1\left[\mu _T\right]\left(x\right)\Longrightarrow L_1\left[\mu
  \right]\left(x\right)=u\left(x\right)
$$
 and the same limit has the statistic $ \xi _T(\hat\vartheta
_T,x)$. Therefore it is sufficient to find such transformation 
$L_2\left[\xi _T(\hat\vartheta _T,\cdot )\right]\left(x\right) $
 that its limit 
is a  Wiener process, say, $L_2\left([U\left(\cdot
  \right)\right]\left(t\right)=w_t $. Below we omit $\vartheta $ in
$h\left(\vartheta ,t\right)$ and denoted the matrix
$$
\NN\left(t\right)=\int_{t}^{1}h\left(\vartheta ,s\right)h^*\left(\vartheta
,s\right){\rm d}s =N\left(\vartheta ,F^{-1}\left(\vartheta ,t\right)\right).
$$
This transformation is given in the following theorem. 
\begin{theorem}
\label{T2}
Suppose that $h\left(s\right)$ is continuous vector-function and  the matrix
$\NN\left(t\right)$ is nondenerate then 
\begin{equation}
\label{19}
L_2\left([U\left(\cdot \right)\right]\left(t\right)\equiv
U\left(t\right)+\int_{0}^{t} \int_{0}^{s}
{h^*\left(v\right)}\;\NN\left(t\right)^{-1}\; h\left(s\right) \,{\rm
  d}U\left(v\right)\;{\rm d}s= w_t
\end{equation}
\end{theorem}
{\bf Proof.} The proof  will be done in several steps. 

{\it Step 1.} We itroduce a Gaussian  process
\begin{equation}
\label{20}
M_t=\int_{0}^{t}q\left(t,s\right)\,{\rm d}U\left(s\right),\qquad 0\leq t\leq 1,
\end{equation}
where the function $q\left(t,s\right)$ we choose as solution of special Fredholm
equation.

{\it Step 2.}  Then we show  that with such choice of $q\left(t,s\right)$ the
process $M_t$
becames a martingale and admits the representation 
$$
M_t=\int_{0}^{t}q\left(s,s\right)\,{\rm d}w_s,\qquad 0\leq t\leq 1,
$$ 
where $w_s,0\leq s\leq 1$ is some Wiener process.

{\it Step 3.}  This representation allows us to obtain the Wiener process
by inverting this equation
$$
w_t=\int_{0}^{t}\frac{1}{q\left(s,s\right)}\,{\rm d}M_s=
U\left(t\right)+\int_{0}^{t}\frac{1}{q\left(s,s\right)}\int_{0}^{s}q'_s\left(s,v\right)
\,{\rm
  d}U\left(v\right)\,{\rm
  d}s,\quad 0\leq t\leq 1.  
$$
This last equality provides us the linear transformation
$$
L_2\left[U\right]\left(t\right)=U\left(t\right)+\int_{0}^{t}\frac{1}{q\left(s,s\right)}\int_{0}^{s}q'_s 
\left(s,v\right)
\,{\rm d}U\left(v\right)\,{\rm d}s =w_t,
$$
and we show that it is equivalent to \eqref{19}.

Now we realize this program. Suppose that $q\left(t,s\right)$ is some
continuous function and the process $M_t$ is defined by the equality
\eqref{20}. Then the  
correlation function of $M_t$ is ($s<t$)
\begin{align*}
R\left(t,s\right)&=\Ex \left[M_tM_s\right]=\Ex
\left[\int_{0}^{t}q\left(t,u\right)\,{\rm
    d}w\left(u\right)-\int_{0}^{t}{q\left(t,u\right)\,\langle \zeta _*,
    h\left(u\right)\rangle }\,{\rm d} u\right]\\
 &\quad
\left[\int_{0}^{s}q\left(s,v\right)\,{\rm
    d}w\left(v\right)-\int_{0}^{s}{q\left(s,v\right)\,\langle \zeta _*, 
    h\left(v\right)\rangle }\,{\rm d}  
  v\right]\\ 
&=\int_{0}^{s}q\left(t,u\right)q\left(s,u\right)\,{\rm
  d}u- \langle\int_{0}^{s}{q\left(s,v\right)\,h\left(v\right)}\,{\rm d}
v,\int_{0}^{t}{q\left(t,u\right)\,h\left(u\right)}\,{\rm d}
u\rangle\\ 
&=\int_{0}^{s}q\left(s,u\right)\left[q\left(t,u\right)-{}\;
  \int_{0}^{t}{q\left(t,v\right)\,\langle
    h\left(u\right),h\left(v\right)}\rangle {\rm d} v\right]\,{\rm d}u. 
\end{align*}

Therefore, if we take $q\left(t,s\right)$ such that it solves  the Fredholm
equation ($t$ is fixed)
\begin{align}
\label{F}
q\left(t,s\right)-{}
\int_{0}^{t}{q\left(t,v\right)\,\langle
  h\left(s\right),h\left(v\right)}\rangle\;{\rm     d} v=1,\qquad
s\in\left[0,t\right], 
\end{align}
then
\begin{align}
\label{mart}
\Ex \left[M_tM_s\right]=\Ex \left[M_s^2\right]=\int_{0}^{s}q\left(s,u\right)\,{\rm d}u.
\end{align}

The solution $q\left(t,s\right) $ of the equation \eqref{F} can be found as
follows. Let us put  
$$
q\left(t,s\right)=1+\langle \int_{0}^{t}q\left(t,v\right)h\left(v\right)\,{\rm
  d}v,{h\left(s\right)}\rangle=1+\langle
A\left(t\right),{h\left(s\right)}\rangle=1+
{h\left(s\right)^*}A\left(t\right), 
$$
where the vector-function $A\left(t\right)$ itself is solution of the following
equation (after multilying \eqref{F} by
${h\left(s\right)}$ and integrating) 
\begin{align*}
A\left(t\right)-\int_{0}^{t}{ h\left(s\right)
}\,h\left(s\right)^*\;{\rm d}s\; A\left(t\right)=\int_{0}^{t}{h\left(s\right)}\;{\rm d}s.
\end{align*}
We can write 
$$
\left({I_d-\int_{0}^{t} {h\left(s\right)}h\left(s\right)^*\;{\rm
    d}s}\right)\; A\left(t\right)= \NN\left(t\right) A\left(t\right)
=\int_{0}^{t}{h\left(s\right)}\;{\rm d}s 
$$ 
($I_d$ is $d\times d$ identity matrix) and remind that $\NN\left(t\right) $
is nondegenerate, then we obtain
\begin{align*}
A\left(t\right)= \NN\left(t\right)^{-1}{\int_{0}^{t}{h\left(s\right)}\;{\rm
  d}s}.
\end{align*}
Therefore, the solution of \eqref{F} is the function 
\begin{align}
\label{sol}
q\left(t,s\right)=1+\langle \NN\left(t\right)^{-1}{\int_{0}^{t}{h\left(v\right)}\;{\rm
  d}v}, h\left(s\right)\rangle  .
\end{align}

The last integral in \eqref{mart} has the following property. 

\begin{lemma}
\label{L}
\begin{equation}
\label{eq}
\int_{0}^{t}q\left(t,s\right)\,{\rm d}s=\int_{0}^{t}q\left(s,s\right)^2{\rm d}s.
\end{equation}
\end{lemma}
{\bf Proof.} 
We  show that 
$$
\frac{{\rm d}}{{\rm d}t}\int_{0}^{t}q\left(t,s\right){\rm d}s=\frac{{\rm
    d}}{{\rm d}t}\int_{0}^{t}q\left(s,s\right)^2{\rm d}s=q\left(t,t\right)^2. 
$$
We have
\begin{align*}
&\frac{{\rm d}}{{\rm d}t}\int_{0}^{t}q\left(t,s\right){\rm d}s=1+\frac{{\rm
    d}}{{\rm d}t}
\left[\int_{0}^{t}h^*\left(s\right){\rm d}s\;\NN\left(t\right)^{-1}
\int_{0}^{t}h\left(v\right){\rm d}v\right]  \\
& \qquad =1+2h^*\left(t\right) \;\NN\left(t\right)^{-1}
\int_{0}^{t}h\left(v\right){\rm d}v\\
&\qquad \quad   + 
\int_{0}^{t}h^*\left(s\right){\rm
  d}s \; \NN\left(t\right)^{-1}  h\left(t\right) h^*\left(t\right)
N\left(t\right)^{-1} \int_{0}^{t}h\left(v\right){\rm d}v\\ 
&\qquad
=\left[1+h^*\left(t\right)\;\NN\left(t\right)^{-1}\int_{0}^{t}h\left(s\right){\rm
    d}s\; 
\right]^2 =q\left(t,t\right)^2.
\end{align*}

The next step is the following Lemma. 
\begin{lemma}
\label{L43-9}
If the Gaussian process $M_s$  satisfies \eqref{mart} and we
have \eqref{eq}  with some continuous  positive function $q\left(s,s\right)$, then 
\begin{align*}
z\left(t\right) = \int_{0}^{t} q\left(s,s\right)^{-1}{\rm d}M_s
\end{align*}
is a Wiener process. 
\end{lemma}
{\bf Proof.}  Consider the partition
$0=s_0<s_1<\ldots<s_N=1$ and put
$$
z_N\left(t\right)=\sum_{s_l\leq t}^{}
q\left(s_{l-1},s_{l-1}\right)^{-1}\left[M_{s_{l}}-M_{s_{l-1}}\right] .
$$
Note that by \eqref{mart}  we have $\Ex M_{s}M_{t} =\Ex M_{s}^2$  for
$s<t$. Hence  for  $l\not=m$
$$
\Ex
\left[M_{s_{l}}-M_{s_{l-1}}\right]\left[M_{s_{m}}-M_{s_{m-1}} \right]
=0. 
$$
This allows us to write 
\begin{align*}
\Ex z_N\left(t\right)z_N\left(s\right)&= \sum_{s_l\leq s}^{}
q\left(s_{l-1},s_{l-1}\right)^{-2} \Ex
\left[M_{s_{l}}-M_{s_{l-1}} \right]^2\\
& = \sum_{s_l\leq s}^{}
q\left(s_{l-1},s_{l-1}\right)^{-2}   \Ex \left[M_{s_{l}}^2-M_{s_{l-1}}^2\right]\\
& = \sum_{s_l\leq s}^{}
q\left(s_{l-1},s_{l-1}\right)^{-2} \int_{s_{l-1}}^{s_l}q\left(v,v\right)^{2}{\rm d}v\longrightarrow s
\end{align*}
as $\max \left|s_l-s_{l-1}\right|\rightarrow 0$.  The same time
$z_N\left(t\right)\rightarrow z\left(t\right)$ in mean-square. Therefore,  $\Ex
z\left(t\right)=0$, $\Ex
z\left(t\right)z\left(s\right)=t\wedge s$ and $z\left(t\right)$ is a Wiener
process $w_t$.

Hence
$$
M_t=\int_{0}^{t}q\left(s,s\right)\,{\rm d}w_s,\qquad t\in [0,1)
$$
is a Gaussian martingale. 
This implies the equality
\begin{align*}
w_t&=\int_{0}^{t}\frac{1}{q\left(s,s\right)}\;{\rm
  d}M_s=U\left(t\right)+\int_{0}^{t}\frac{1}{q\left(s,s\right)}\;\int_{0}^{s}
q'_s\left(s,v\right) \,{\rm d}U\left(v\right)\;{\rm d}s .
\end{align*}
For the derivative  $q'_t\left(t,s\right) $ we can write 
\begin{align*}
&q'_t\left(t,s\right)= 
\left(A'\left(t\right),h\left(s\right)\right)\\
&\qquad =h^* \left(s\right)   \NN\left(t\right)^{-1}\;h\left(t\right) h^*\left(t\right) 
\NN\left(t\right)^{-1}\;\int_{0}^{t}h\left(v\right){\rm d}v
+h^*\left(s\right)\; \NN\left(t\right)^{-1}\; h\left(t\right)\\
&\qquad =h^*\left(s\right)    \NN\left(t\right)^{-1}\;h\left(t\right)  \left[ h^*\left(t\right)
\NN\left(t\right)^{-1}\;\int_{0}^{t}h\left(v\right){\rm d}v+1\right]\\
&\qquad =h^*\left(s\right)   \NN\left(t\right)^{-1}\;h\left(t\right) \,q\left(t,t\right).
\end{align*}
Hence
\begin{align*}
\frac{q'_s\left(s,v\right)}{q\left(s,s\right)}=h^*\left(v\right)
\NN\left(s\right)^{-1}\;h\left(s\right) 
\end{align*}
and we obtain the final expression 
\begin{align*}
w_t=
U\left(t\right)+\int_{0}^{t} \int_{0}^{s}  {h^*\left(v\right)}\;\NN\left(t\right)^{-1}\;
h\left(s\right) \,{\rm d}U\left(v\right)\;{\rm d}s.
\end{align*}

This is the explicit linear transformation $w_t =
L_2\left[U \right]\left(t\right)$ of the process $U\left(\cdot \right) 
$ in the Wiener 
process $w_t$ and this proves the Theorem \ref{T2}.

Let us  denote 
\begin{align*}
g\left(\vartheta ,y\right)= \frac{\dot S\left(\vartheta
  ,y\right)}{\;\;\sigma \left(y\right)}, \quad   
\NN\left(\vartheta, x \right)= \int_{x}^{\infty
}\frac{\dot S\left(\vartheta 
  ,z\right)\dot S\left(\vartheta 
  ,z\right)^*}{\sigma
  \left(z\right)^2}\,f\left(\vartheta ,z\right)\,{\rm d}z, 
\end{align*}
then we can write 
\begin{align*}
&w_{F\left(\vartheta, x\right)}=
U\left(F\left(\vartheta ,x\right)\right)\\
&\qquad \quad +\int_{-\infty }^{x} \int_{-\infty }^{y} g^*\left(\vartheta
,y\right) {\NN\left(\vartheta ,x 
  \right)}^{-1}{}\; 
g\left(\vartheta ,z\right) \,{\rm d} U\left(F\left(\vartheta
,z\right)\right)\; f\left(\vartheta ,y\right){\rm d}y, 
\end{align*}
i.e., this transformation of $U\left(\cdot \right)$ does not depend on information matrix ${\rm
      I}\left(\vartheta \right) $. Of course, $U\left(\cdot \right)$  itself
depends on ${\rm
      I}\left(\vartheta \right) $.

To construct the test we have to replace $U\left(F\left(\vartheta
      ,x\right)  \right),g\left(\vartheta,y
\right)$ and $\NN\left(\vartheta ,y \right)$   in \eqref{19} by their empirical versions
based on observations only 
$$
\xi _T\left(\hat\vartheta _T,x\right),\quad g\left(\hat\vartheta _T,y
\right)=\frac{\dot S\left(\hat\vartheta _T
  ,y\right)}{\;\;\sigma \left(y\right)},\qquad \NN\left(\hat\vartheta _T ,x \right)
$$  
respectively and to study
\begin{align*}
&v_T\left(\hat\vartheta _T, x\right)=
\xi _T\left(\hat\vartheta _T ,x\right)\\
& \quad +\int_{-\infty }^{x} \int_{-\infty }^{y} g^*\left(\hat\vartheta
_T,y\right) {\NN\left(\hat\vartheta _T ,x
  \right)}^{-1}{}\; 
g\left(\hat\vartheta _T,z\right) \,{\rm d} \xi _T\left(\hat\vartheta _T
,z\right)\; {\rm d}F\left(\hat\vartheta _T ,y\right). 
\end{align*}
Then we have to show that 
$$
v_T(\hat\vartheta _T, x)-v_T\left(\vartheta, x\right)\rightarrow 0,\quad
v_T\left(\vartheta, x\right)\Longrightarrow  w_{F\left(\vartheta, x\right)}.
$$
Unfortunately we can not do it directly.  We have to avoid the calculation of the integral
$$
S\left(\hat\vartheta _T,y\right)=\int_{-\infty }^{y}
     g\left(\hat\vartheta _T,z
\right) \,{\rm d}\xi _T\left(\hat\vartheta _T,z\right)
$$
because this integral is equivalent in some sense to the It\^o stochastic
integral and $\hat\vartheta _T$ depends on the whole trajectory $\left(X_t,0\leq
t\leq T\right)$. One way is to use the discrete approximation of this integral
\begin{align*}
K_n\left(\hat\vartheta _T,y\right)=\sum_{z_i<y}^{}g\left(\hat\vartheta _T,z_i
\right) \,\left[\xi _T\left(\hat\vartheta _T,z_{i+1}\right) -\xi
  _T\left(\hat\vartheta _T,z_{i}\right)\right] 
\end{align*}
and to show that 
$$
K_n\left(\hat\vartheta _T,y\right)-K_n\left(\vartheta ,y\right)\rightarrow
0,\qquad K_n\left(\vartheta ,y\right)-K\left(\vartheta ,y\right)\rightarrow 0. 
$$
Another possibility is to replace the corresponding stochastic integral by the
ordinary one what we do below.

Introduce two functions 
\begin{align*}
Q\left(\vartheta ,x,y\right)&= \int_{y\wedge x }^{x} g\left(\vartheta
,v\right) {\NN\left(\vartheta ,x \right)}^{-1}{}\;{\rm d}F\left(\vartheta
,v\right),\\ 
R\left(\vartheta ,x,y\right)&= \frac{\langle \dot
  S\left(\vartheta ,y\right),Q\left(\vartheta, x,y\right)\rangle}{\sigma
  \left(y\right)^2}
\end{align*}
and the statistic
\begin{align*}
V_T\left(\hat\vartheta _T,x\right)&=\xi _T\left(\hat\vartheta_T,x\right)
-\frac{1}{2\sqrt{T}}\int_{0}^{T} 
  \left[R_y'\left(\hat\vartheta _T ,x,X_s\right)\sigma 
\left(X_s\right)^2{\rm   d}s\right.\\
&\qquad \left.+2{R\left(\hat\vartheta _T
  ,x,X_s\right)}\,S(\hat\vartheta _T,X_s)\right]\,{\rm d}s.
\end{align*}

The main result of this work is the following theorem. 
\begin{theorem}
\label{T3} 
Let the conditions ${\cal ES}, {\cal A}_0$ and ${\cal R}_1,{\cal R}_2$ be
fulfilled, then the test {\rm $\hat\psi _T=\1_{\left\{\delta _T>c_\varepsilon
    \right\}}$} with $\delta  _T=$ and $c_\varepsilon $ defined by the relations 
\begin{align}
\label{test}
\delta  _T=\int_{-\infty }^{\infty }V_T\left( \hat\vartheta _T,x\right)^2{\rm
  d}F(\hat \vartheta _T,x),\qquad \Pb\left(\int_{0}^{1}w_t^2{\rm
  d}t>c_\varepsilon \right)=\varepsilon  
\end{align}
is ADF and belongs to ${\cal K}_\varepsilon $.
\end{theorem}
{\bf Proof.}

Let us suppose that $m\left(\vartheta ,z \right)$ is piece-wise continuous
function and consider the calculation of the integral
$$
\int_{a }^{b}
     g\left(\vartheta ,z
\right) \,{\rm d}\xi _T\left(\vartheta,z\right).
$$
  For any partition $a=z_1<z_2\ldots < z_K=b$ and
$\max\left|z_{k+1}-z_k\right|\rightarrow 0$ we have
\begin{align*}
&\sum_{k=1}^{K-1}g\left(\vartheta ,\tilde z_k\right) \left[
    \xi _T\left(\vartheta,z_{k+1}\right)-\xi _T\left(\vartheta,z_{k}\right)\right]\\ 
&\qquad
  =\frac{1}{\sqrt{T}}\int_{0}^{T} \frac{\sum_{k=1}^{N-1}g\left(\vartheta
    ,\tilde z_k\right)\1_{\left\{z_k\leq  X_s<z_{k+1}\right\}}}{\sigma
    \left(X_s\right)}\,{\rm d}X_s\\ 
&\qquad\qquad \quad -\frac{1}{\sqrt{T}}\int_{0}^{T} \frac{\sum_{k=1}^{N-1}g\left(\vartheta
    ,\tilde z_k\right)S(\vartheta,X_s)\1_{\left\{z_k\leq  X_s< z_{k+1}\right\}}}{\sigma
    \left(X_s\right)}\,{\rm d}s
  \\
 &\qquad\longrightarrow \frac{1}{\sqrt{T}}\int_{0}^{T} \frac{g\left(\vartheta
    ,X_s\right)\1_{\left\{a\leq  X_s< b\right\}}}{\sigma
    \left(X_s\right)}\,{\rm d}X_s\\ 
&\qquad\qquad \quad -\frac{1}{\sqrt{T}}\int_{0}^{T} \frac{g\left(\vartheta
    ,X_s\right)S(\vartheta,X_s)\1_{\left\{a\leq  X_s< b\right\}}}{\sigma
    \left(X_s\right)}\,{\rm d}s
\end{align*}
Therefore we have the equality
\begin{align}
  &\int_{-\infty }^{y}\frac{\dot S\left(\vartheta ,z\right)}{\sigma (z)}\,{\rm
    d}\xi _T\left(\vartheta ,z\right) 
=\frac{1}{\sqrt{T}}\int_{0}^{T} \frac{\dot S\left(\vartheta
  ,X_s\right)\1_{\left\{ X_s< y\right\}}}{\sigma 
    \left(X_s\right)^2}\,{\rm d}X_s\nonumber\\ 
&\qquad  \quad -\frac{1}{\sqrt{T}}\int_{0}^{T} \frac{\dot S\left(\vartheta
    ,X_s\right)S(\vartheta,X_s)\1_{\left\{ X_s< y\right\}}}{\sigma
    \left(X_s\right)^2}\,{\rm d}s.
\label{77}
\end{align}
Further, by Fubini theorem
\begin{align*}
&J_T\left(\vartheta ,x\right)=\int_{-\infty }^{x}  g^*\left(\vartheta
  ,y\right) {\NN\left(\vartheta ,x 
  \right)}^{-1}{}\; \int_{-\infty }^{y}
g\left(\vartheta ,z\right) \,{\rm d} \xi _T\left(\vartheta,z\right)\;
{\rm d}F\left(\vartheta ,y\right),\\
&\quad =\frac{1}{\sqrt{T}}\int_{0}^{T}  \frac{\dot S\left(\vartheta ,X_s\right)^*}{\sigma
    \left(X_s\right)^2}{\NN\left(\vartheta ,x
  \right)}^{-1} \int_{X_s\wedge x }^{x} g\left(\vartheta ,y\right)   {\rm
  d}F\left(\vartheta ,y\right)\,{\rm d}X_s\\ 
&\qquad -\frac{1}{\sqrt{T}}\int_{0}^{T}\frac{\dot S\left(\vartheta
    ,X_s\right)^*S(\vartheta,X_s)}{\sigma
    \left(X_s\right)^2} {\NN\left(\vartheta ,x
  \right)}^{-1} \int_{X_s\wedge x }^{x}  g\left(\vartheta ,y\right) {}\;{\rm
  d}F\left(\vartheta ,y\right)\;{\rm d}s\\ 
&\quad =\frac{1}{\sqrt{T}}\int_{0}^{T} R\left(\vartheta ,x,X_s\right)\,{\rm
  d}X_s-\frac{1}{\sqrt{T}}\int_{0}^{T} {R\left(\vartheta
  ,x,X_s\right)}\,S(\vartheta,X_s)\,{\rm d}s. 
\end{align*}

By It\^o formula 
\begin{align*}
\int_{0}^{T} R\left(\vartheta  ,x,X_s\right)\,{\rm
  d}X_s=\int_{X_0}^{X_T} R\left(\vartheta ,x,y\right)\,{\rm
  d}y-\frac{1}{2}\int_{0}^{T} R_y'\left(\vartheta ,x,X_s\right)\sigma
\left(X_s\right)^2{\rm   d}s.
\end{align*}
Hence we have no more stochastic integrals and can substitute the estimator 
\begin{align*}
&\sqrt{T}J_T\left(\hat\vartheta_T ,x\right)=\int_{X_0}^{X_T}
  R\left(\hat\vartheta_T ,x,y\right)\,{\rm d}y\\
&\quad
-\int_{0}^{T} \left[{R\left(\hat\vartheta_T 
  ,x,X_s\right)}S(\hat\vartheta_T ,X_s)+\frac{1}{2}R_y'\left(\hat\vartheta_T  ,x,X_s\right)\sigma
\left(X_s\right)^2 \right]{\rm d}s\\
&\quad =\int_{X_0}^{X_T}
  R\left(\hat\vartheta_T ,x,y\right)\,{\rm d}y\\
&\quad
-\int_{0}^{T} \left[{R\left(\hat\vartheta_T 
  ,x,X_s\right)}S(\vartheta ,X_s)+\frac{1}{2}R_y'\left(\hat\vartheta_T  ,x,X_s\right)\sigma
\left(X_s\right)^2 \right]{\rm d}s\\
&\quad
+\int_{0}^{T} {R\left(\hat\vartheta_T 
  ,x,X_s\right)}\left[S(\vartheta ,X_s) -S(\hat\vartheta_T ,X_s)\right]{\rm d}s.
\end{align*}
 Further (below $\hat u_T=\sqrt{T}\left(\hat\vartheta_T-\vartheta\right) $)
\begin{align*}
&\left[J_T\left(\hat\vartheta_T ,x\right)-J_T\left(\vartheta
    ,x\right)\right]=\langle\frac{\hat u_T}{{T}},\int_{X_0}^{X_T} 
  \dot R\left(\vartheta,x,y\right)\,{\rm d}y\rangle\\
&\quad
-\langle\frac{\hat u_T}{{T}},\int_{0}^{T} \left[{\dot R\left(\vartheta
  ,x,X_s\right)}S(\vartheta ,X_s)+\frac{1}{2}\dot R_y'\left(\vartheta  ,x,X_s\right)\sigma
\left(X_s\right)^2 \right]{\rm d}s\rangle\\
&\quad
-\langle\frac{ \hat u_T}{{T}},\int_{0}^{T} {R\left(\vartheta 
  ,x,X_s\right)}\dot S(\vartheta ,X_s) {\rm d}s\rangle+o\left(1\right).
\end{align*}
Note that by Theorem 2.8 \cite{Kut04} for any $p>0$
\begin{equation}
\label{66}
\sup_\vartheta \Ex_\vartheta \left|\hat\vartheta _T-\vartheta \right|^p\leq C\,T^{-\frac{p}{2}}.
\end{equation}
Using once more the It\^o formula we obtain
\begin{align*}
&\int_{X_0}^{X_T} 
  \dot R\left(\vartheta,x,y\right)\,{\rm d}y-\int_{0}^{T} \left[{\dot R\left(\vartheta
  ,x,X_s\right)}S(\vartheta ,X_s)+\frac{1}{2}\dot R_y'\left(\vartheta  ,x,X_s\right)\sigma
\left(X_s\right)^2 \right]{\rm d}s\\
&\qquad =\int_{0}^{T} {\dot R\left(\vartheta
  ,x,X_s\right)}{\rm d}W_s.
\end{align*}
Hence
\begin{align*}
&\left(\Ex_\vartheta \langle\frac{\hat u_T}{T},\int_{0}^{T} {\dot R\left(\vartheta
  ,x,X_s\right)}\;{\rm d}W_s\rangle\right)^2\\
&\qquad \leq \Ex_\vartheta \left|\hat
u_T\right|^2\;\left| \frac{1}{T}\int_{0}^{T} {\dot R\left(\vartheta
  ,x,X_s\right)}\;{\rm d}W_s  \right|^2\leq \frac{C}{{T}}
 ,
\end{align*}
and we can write
\begin{align*}
J_T\left(\hat\vartheta_T ,x\right)&=\frac{1}{\sqrt{T}}\int_{0}^{T} { R\left(\vartheta
  ,x,X_s\right)}\;{\rm d}W_s\\
&\quad -\frac{ 1}{{T}}\int_{0}^{T} {R\left(\vartheta 
  ,x,X_s\right)}\langle\hat u_T, \dot S(\vartheta ,X_s)\rangle\; {\rm d}s+o\left(1\right)
\end{align*}
Therefore
\begin{align*}
&V_T\left(\hat\vartheta _T,x\right)= \xi _T\left( \hat\vartheta
_T,x\right)+\frac{1}{\sqrt{T}}\int_{0}^{T} { R\left(\vartheta 
  ,x,X_s\right)}\;{\rm d}W_s\\
&\quad -\frac{ 1}{{T}}\int_{0}^{T} {R\left(\vartheta 
  ,x,X_s\right)} \langle \hat u_T,\dot S(\vartheta ,X_s)\rangle\; {\rm d}s+o\left(1\right)=\hat V_T\left(\hat\vartheta _T,x\right)+o\left(1\right),
\end{align*}
where we put
\begin{align*}
\hat V_T\left(\hat\vartheta _T,x\right)&=\xi _T\left( \hat\vartheta
_T,x\right)+\frac{1}{\sqrt{T}}\int_{0}^{T} { R\left(\vartheta
  ,x,X_s\right)}\;{\rm d}W_s\\ &\quad -\frac{ 1}{{T}}\int_{0}^{T}
     {R\left(\vartheta ,x,X_s\right)}\langle \hat u_T ,\dot S(\vartheta
     ,X_s)\rangle\; {\rm d}s. 
\end{align*}
To prove the convergence 
\begin{align*}
\delta _T&=\int_{-\infty }^{\infty }\hat V_T\left(\hat\vartheta
_T,x\right)^2{\rm d}F\left(\hat\vartheta _T,x\right)+o\left(1\right)\\
&\qquad 
\Longrightarrow \int_{-\infty }^{\infty }w_{F\left(\vartheta ,x\right)} ^2{\rm
  d}F\left(\vartheta,x\right)=\int_{0}^{1}w_t^2{\rm d}t 
\end{align*}
we have to verify the following properties:
\begin{enumerate}
\item  For any $x_1,\ldots,x_k$
$$
\left(\hat V_T(\hat\vartheta _T,x_1),\ldots,\hat V_T(\hat\vartheta
_T,x_k)\right)\Longrightarrow \left(w_{F\left(\vartheta
  ,x_1\right)},\ldots,w_{F\left(\vartheta ,x_k\right)}\right). 
$$
\item For any $\delta >0$ there exist $L>0$ such that
\begin{equation}
\label{22}
\int_{\left|x\right|>L}^{}\Ex_\vartheta\hat V_T(\hat\vartheta
_T,x)^2f(\hat\vartheta _T,x)\,{\rm d}x<\delta  .
\end{equation}
\item For $\left|x_i\right|<L,i=1,2$
\begin{equation}
\label{23}
\Ex_\vartheta\left|\hat  V_T(\hat\vartheta
_T,x_2)-\hat  V_T(\hat\vartheta
_T,x_1)\right|^2\leq C\;\left|x_2-x_1\right|^{1/2}.
\end{equation}
\end{enumerate}

The first convergence follows from \eqref{55}, central limit theorem for
stochastic integrals and the law of large numbers 
$$
\frac{ 1}{T}\int_{0}^{T}
     R\left(\vartheta ,x_i,X_s\right) \dot S(\vartheta ,X_s)\; {\rm
       d}s\longrightarrow 
\int_{-\infty }^{\infty }R\left(\vartheta ,x_i,y\right)\dot S(\vartheta
,y)f\left(\vartheta ,y\right){\rm d}y.
$$
Here $i=1,\ldots,k. $
Indeed, we obtain the joint asymptotic normality 
\begin{align*}
\hat V_T(\hat\vartheta_T,x_i)&\Longrightarrow u\left(x\right)+\int_{-\infty
}^{\infty }R\left(\vartheta ,x_i,y\right){\rm
  d}W\left(F\left(\vartheta ,y\right)  \right)\\
&\qquad -\int_{-\infty }^{\infty }R\left(\vartheta ,x_i,y\right)
    \langle {\rm I}\left(\vartheta \right)^{-1}\Delta \left(\vartheta \right),
    \dot S(\vartheta 
,y)\rangle {\rm d}F\left(\vartheta ,y\right). 
\end{align*}
Note that the limit of \eqref{77} is equivalent to 
\begin{align*}
\int_{-\infty }^{y}\frac{\dot S\left(\vartheta ,z\right)}{\sigma
  \left(z\right)}{\rm d}u\left(x\right)&=\int_{-\infty }^{y}\frac{\dot
  S\left(\vartheta ,z\right)}{\sigma 
  \left(z\right)}{\rm d}W\left(F\left(\vartheta ,z\right)\right)\\
&\qquad -  \int_{-\infty }^{y}\langle {\rm
  I}\left(\vartheta \right)^{-1}\Delta \left(\vartheta \right),\dot
S\left(\vartheta ,z\right)\rangle\frac{\dot S\left(\vartheta ,z\right)}{\sigma 
  \left(z\right)^2} {\rm d}F\left(\vartheta ,z\right).
\end{align*} 

To check \eqref{22} we write
\begin{align*}
\Ex_\vartheta \xi _T\left(x\right)^2&\leq 2 \Ex_\vartheta
\left(\frac{1}{\sqrt{T}}\int_{0}^{T}{\1_{\left\{X_s<x\right\}} }{}{\rm
  d}W_s\right)^2\\
&\qquad +2\Ex_\vartheta\left( \langle\hat
u_T\;,\frac{1}{T}\int_{0}^{T}\frac{\1_{\left\{X_s<x\right\}} \dot
  S\left(\tilde \vartheta _T,X_s\right) }{\sigma \left(X_s\right)}{\rm
  d}s\rangle \right)^2\\ 
&\leq 2F\left(\vartheta ,x\right)+2\Ex_\vartheta\left|\hat u_T\right|^2
\left|\frac{1}{T}\int_{0}^{T}\frac{\1_{\left\{X_s<x\right\}} \dot
  S\left(\tilde \vartheta _T,X_s\right) }{\sigma \left(X_s\right)}{\rm
  d}s\right|^2\leq C. 
\end{align*}
Remind that  by conditions ${\cal A}_0,{\cal R}_1,{\cal R}_2,$  all
related functions have
polynomial  majorants. The invariant density $f\left(\vartheta ,x\right)$ by
condition ${\cal A}_0$ has exponentially decreasing tails: there exist the
constants $c_1>0,C_2>0$ such that
$$
f\left(\vartheta ,x\right)\leq C_2\,e^{-c_2\left|x\right|}.
$$
 Therefore all 
mathematical expectations are finite. 

Further, 
\begin{align*}
&\Ex_\vartheta \left|\hat V_T\left(\hat\vartheta _T,x_2\right)-\hat
V_T\left(\hat\vartheta _T,x_1\right)\right|^2 \leq 3\Ex_\vartheta \left|\xi
_T\left(x_2\right)-\xi _T\left(x_1\right)\right|^2 \\
&\qquad
+3\Ex_\vartheta\left|\frac{1}{\sqrt{T}}\int_{0}^{T}\left[R\left(\vartheta
  ,x_2,X_s\right)-R\left(\vartheta ,x_1,X_s\right)\right] {\rm d}W_s \right|^2
\\ 
&\qquad +3\Ex_\vartheta\left|\frac{1}{{T}}\int_{0}^{T}\left[R\left(\vartheta
  ,x_2,X_s\right)-R\left(\vartheta ,x_1,X_s\right)\right] \langle \hat u_T, \dot
  S\left(\tilde \vartheta _T,X_s\right)\rangle\,  {\rm d}s  \right|^2\\
&\quad \leq C\left(L\right)\,\left|x_2-x_1\right|^{1/2}.
\end{align*}
For example, ($x_1<x_2$)
\begin{align*}
&\Ex_\vartheta \left|\xi
_T\left(x_2\right)-\xi _T\left(x_1\right)\right|^2\leq 2 \Ex_\vartheta
\left(\frac{1}{\sqrt{T}}\int_{0}^{T}{\1_{\left\{x_1<X_s<x_2\right\}}}{}{\rm
  d}W_s\right)^2\\
&\qquad +2\Ex_\vartheta\left( \langle\hat
u_T\;,\frac{1}{T}\int_{0}^{T}\frac{\1_{\left\{x_1<X_s<x\right\}} \dot
  S\left(\tilde \vartheta _T,X_s\right) }{\sigma \left(X_s\right)}{\rm
  d}s\rangle \right)^2\\ 
&\qquad \leq 2\int_{x_1}^{x_2}f\left(\vartheta ,y\right)\,{\rm d}y +2
\left(\Ex_\vartheta\left|\hat u_T\right|^4\right)^{1/2} \left(
\int_{x_1}^{x_2}P\left(y\right)f\left(\vartheta ,y\right)\,{\rm d}y   \right) ^{1/2}\\ 
&\qquad \leq C\left|x_2-x_1\right|+ C\left|x_2-x_1\right|^{1/2}\leq
C\left(L\right)\left|x_2-x_1\right|^{1/2}. 
\end{align*}
Here $P\left(y\right)$ is some polynome. 
These properties of $V_T(\hat\vartheta
_T,x) $ allow us (see Theorem A1.22 \cite{IH81}) to verify the convergence
$$
\int_{-\infty }^{\infty } V_T(\hat\vartheta
_T,x)^2f(\hat\vartheta _T,x)\,{\rm d}x\Longrightarrow \int_{-\infty }^{\infty
} w_{F\left(\vartheta ,x\right)}^2\;{\rm d}F\left(\vartheta
,x\right)=\int_{0}^{1}w_t^2\;{\rm d}t. 
$$

{\bf Example.} {\it Linear case.}
Let us consider the one-dimensional ($d=1$) linear case 
$$
{\rm d}X_s=\vartheta a\left(X_s\right){\rm d}s + \sigma \left(X_s\right){\rm
  d}W_s,\quad X_0,\quad 0\leq s\leq T .
$$ 
We have some simplification because we have no more problem with the
calculation of stochastic integral and the statistic can be calculated as
follows.  Let us denote
\begin{align*}
B_T\left(\hat\vartheta _T,x\right)= \xi _T(\hat\vartheta _T,x) +\int_{-\infty }^{x} \frac{
  a\left(y\right)A_T(\hat\vartheta _T,y)}{\NN\left(\hat\vartheta _T,y\right)\sigma
  \left(y\right)}\; {\rm d}F(\hat\vartheta _T,y),
\end{align*}
where
$$
 \NN\left(\vartheta 
,y\right)=\int_{y}^{\infty }\frac{a\left(z\right)^2}{\sigma
  \left(z\right)^2}\,f\left(\vartheta ,z\right)\;{\rm d}z 
$$
and (see \eqref{77})
\begin{align*}
A_T(\hat\vartheta _T,y)& = \frac{1}{\sqrt{T}}\int_{0}^{T}
\frac{a\left(X_s\right)\1_{\left\{ X_s< y\right\}}}{\sigma
  \left(X_s\right)^2}\,\left[{\rm d}X_s -\hat\vartheta _Ta\left(X_s\right)\,{\rm d}s  \right].
\end{align*}
Then we obtain the convergence 
$$
\delta _T=\int_{-\infty }^{\infty } B_T(\hat\vartheta _T,x)^2{\rm
  d}F(\hat\vartheta _T,x)\Longrightarrow  \int_{0}^{1}w_t^2\;{\rm d}t
$$
Hence the test $\hat\psi _T=\1_{\left\{\delta _T>c_\varepsilon \right\}}$ is  ADF.

\section{Discussion}

 In Theorem \ref{T2} the condition of existence of the finite solution on the interval
$[0,1) $ is the following: for any $t\in \left(0,1\right) $  the matrix
\begin{equation}
\label{sc}
\NN\left(t\right)=\int_{t}^{1} h\left(v\right)h^*\left(v\right){\rm d}v
\end{equation}
is positive defined. 
Of course, we have  to check it for {\it any close to} 1 value of $t <1$. The quantity
$\NN\left(t\right)={\rm I}\left(\vartheta \right)^{-1} {\rm I}_t\left(\vartheta
\right)$, where
$$
{\rm I}_t\left(\vartheta\right)=\int_{x}^{\infty } \frac{\dot
S\left(\vartheta ,y\right)\dot
S\left(\vartheta ,y\right)^*}{\sigma \left(y\right)^2}f\left(\vartheta ,y\right) {\rm d}y
$$
 is the Fisher information in the case of censored
observations
 $$
Y_t=X_t\;\1_{\left\{X_t>x\right\}}
$$
and the condition \eqref{sc} means that this Fisher information is positive defined.

For example, if $d=1$ and we suppose  that 
$$
h\left(1\right)=\lim_{t\rightarrow 1} \frac{\left|\dot
  S\left(\vartheta ,F^{-1}\left(\vartheta ,t\right)\right)\right|}{\sigma
  \left(F^{-1}\left(\vartheta ,t\right)\right)\sqrt{{{\rm 
        I}\left(\vartheta \right)}}} =\lim_{y\rightarrow \infty } \frac{\left|\dot
  S\left(\vartheta ,y\right)\right|}{\sigma \left(y\right)\sqrt{{{\rm
        I}\left(\vartheta \right)}}} >0,
$$
then the condition \eqref{sc} is fulfilled.

It is easy to see that for Ornstein-Uhlenbeck process $h\left(1\right)=\infty
$, but the integral of $h\left(\cdot \right)^2$ on $\left[0,1\right]$ is
finite and equal to 1.

 Note that if the function $\dot S\left(\vartheta ,y\right)=0$ for $y\geq b$
 with some $b$, then we have finite solution $q\left(t,s\right),s\in
 \left[0,t\right]$ for the values $t\in[0,F\left(\vartheta, b\right))$ only.

The transformation $L\left[\cdot \right]$ of the limit process \eqref{19}  coincides with one
proposed in  \cite{Kh81} and the difference is in the 
proofs. The transformation $L\left[\cdot \right]$ in \cite{Kh81} is based on
two strong results: one is due to  Hitsuda \cite{Hit68}, which gives the linear
representation of a Gaussian process with measure equivalent to the measure of
Wiener process and the second is due to  Shepp \cite{Sh}, which gives the
condition of equivalence of the process \eqref{bp} $U\left(\cdot \right)$ on
any interval $\left[0,\tau \right], \tau <1$ to the Wiener process $W_s, 0\leq
t\leq \tau $ and then $\tau \rightarrow 0$.  We do not use these two results
and give the direct martingale 
proof using  the solution of Fredholm equation of the second kind with
degenerated kernel.  

{\bf Aknowledgement.} The authors are deeply gratefull to R. Liptser for
fruitful discussions.

\end{document}